\documentclass[english]{article}

\usepackage[T1]{fontenc}
\usepackage[latin9]{inputenc}
\usepackage{color}
\usepackage{amsmath}
\usepackage{mathrsfs}
\usepackage{babel}
\usepackage[active]{srcltx}
\usepackage{amssymb}
\usepackage{a4wide}

\makeatletter

\providecommand{\U}[1]{\protect\rule{.1in}{.1in}}
\providecommand{\U}[1]{\protect\rule{.1in}{.1in}}
\newtheorem{theorem}{Theorem}

\newtheorem{remark}{Remark}

\newenvironment{proof}[1][Proof]{\noindent\textbf{#1.} }{\ \rule{0.5em}{0.5em}}

\makeatother

\begin{document}
\title{The sum theorem for maximal monotone operators in reflexive Banach
spaces revisited}
\author{M.D. Voisei}
\date{{}}
\maketitle
\begin{abstract}
Our goal is to present a new shorter proof for the maximal monotonicity of
the Minkowski sum of two maximal monotone multi-valued operators defined
in a reflexive Banach space under the classical interiority condition
involving their domains.
\end{abstract}
\textbf{Keywords} maximal monotone operator, Minkowski sum

\strut

\noindent \textbf{Mathematics Subject Classification (2010)} 47H05,
46N10.

\section{Preliminaries}

Recall the following sum rule for maximal monotone operators:

\begin{theorem} \label{RM} (Rockafellar \cite[Theorem\ 1\ (a)]{MR0282272})
Let $X$ be a reflexive Banach space with topological dual $X^{*}$
and let $A,\ B:X\rightrightarrows X^{*}$ be multi-valued maximal
monotone operators from $X$ to $X^{*}$. If $D(A)\cap\operatorname*{int}D(B)\neq\emptyset$
then $A+B$ is maximal monotone. Here $D(T):=\{x\in X\mid T(x)\neq\emptyset\}$
is the domain of $T:X\rightrightarrows X^{*}$ and ``$\operatorname*{int}S$'' denotes the
interior of $S\subset X$. \end{theorem}

The proof of \cite[Theorem\ 1]{MR0282272} relies on the use of the
duality mapping $J$ of $X$ and the (Minty's style) characterization
of maximal monotone operators defined in reflexive Banach spaces. Similar
arguments are used in the presence of an improved qualification constraint
in a second proof of Theorem \ref{RM} (see \cite[Corollary 3.5, p. 286]{MR1398155}).
A third proof of the main theorem involves the exact convolution of
some specially constructed functions based on the Fitzpatrick functions
of $A$ and $B$ (see \cite[Corollary\ 4, p. 1166]{MR2160665}). A
different proof of Theorem \ref{RM} is based on the dual-representability
$A+B$ in the presence of the qualification constraint (see \cite[Remark\ 1,\ p.\ 276]{MR2577332})
and the fact that in a reflexive Banach space dual-representability
is equivalent to maximal monotonicity (see e.g. \cite[Theorem\ 3.1,\ p.\ 2381]{MR1974634}).
All the previously mentioned proofs make use of the duality mapping
$J$ which is characteristic to a normed space. 

Our proof relies on the normal cone, is based on full-range characterizations
of maximal monotone operators with bounded domain, and uses the representability of sums
of representable operators, but, avoids the use of $J$ or the norm.
The following intermediary result, is the
main ingredient of our argument. 

\begin{theorem} Let $X$ be a reflexive Banach space, let $T:X\rightrightarrows X^{*}$
be maximal monotone, and let $C\subset X$ be closed convex and bounded.
If $D(T)\cap\operatorname*{int}C\neq\emptyset$ then $T+N_{C}$ is
maximal monotone. Here $N_{C}$ denotes the normal cone to $C$ and
is defined by $x^{*}\in N_{C}(x)$ if, for every $y\in C$, $x^{*}(y-x)\le0$.
\end{theorem}

Recall that a multi-valued operator $T:X\rightrightarrows X^{*}$
is \emph{monotone} if, for every $x_{1}^{*}\in T(x_{1})$, $x_{2}^{*}\in T(x_{2})$,
$\langle x_{1}-x_{2},x_{1}^{*}-x_{2}^{*}\rangle\ge0$. Here $\langle x,x^{*}\rangle:=c(x,x^{*}):=x^{*}(x)$,
$x\in X$, $x^{*}\in X^{*}$.

An element $z=(x,x^{*})\in X\times X^{*}$ is \emph{monotonically
related} (m.r. for short) \emph{to} $T$ if, for every $(a,a^{*})\in\operatorname*{Graph}T:=\{(u,u^{*})\in X\times X^{*}\mid u\in D(T),\ u^{*}\in T(u)\}$,
$\langle x-a,x^{*}-a^{*}\rangle\ge0$.

An operator $T:X\rightrightarrows X^{*}$ is \emph{maximal monotone}
if every m.r. to $T$ element $z=(x,x^{*})\in X\times X^{*}$ belongs
to $\operatorname*{Graph}T$. 

\section{Proofs of the main result}

\strut

\begin{proof}[Proof of Theorem 2] The operator is representable, which follows
from the facts that $T$, $N_{C}$ are maximal monotone thus representable
and $D(T)\cap\operatorname*{int}C\neq\emptyset$ (see e.g. \cite[Corollary\ 5.6]{MR2453098} or \cite[Theorem 16, p. 818]{MR3958031}). 

We prove that $R(T+N_{C})=X^{*}$ which implies that $T+N_C$ is of NI--type and so it is maximal monotone (see \cite[Theorem 3.4, p. 465]{MR2453098} or \cite[Theorem 1 (ii), (7)]{MR2577332}). 

It suffices to prove that $0\in R(T+N_{C})$
otherwise we replace $T$ by $T-x^{*}$ for an arbitrary $x^{*}\in X^{*}$.

Consider $F(x,x^{*}):=\varphi_{T}(x,x^{*})+g(x,x^{*})$, with $g(x,x^{*}):=\iota_{C}(x)+\sigma_{C}(-x^{*})$,
where 
\begin{equation}
\varphi_{T}(x,x^{*}):=\sup\{\langle x-a,a^{*}\rangle+\langle a,x^{*}\rangle\mid(a,a^{*})\in\operatorname*{Graph}T\},\ (x,x^{*})\in X\times X^{*},\label{ff}
\end{equation}
is the Fitzpatrick function of $T$, $\iota_{C}(x)=0$, for $x\in C$;
$\iota_{C}(x)=+\infty$, otherwise, and $\sigma_{C}(x^{*}):=\sup_{x\in C}\langle x,x^{*}\rangle$,
$x^{*}\in X^{*}$.

Then $F\ge0$ due to $\varphi_{T}(x,x^{*})\ge\langle x,x^{*}\rangle$
and $\iota_{C}(x)+\sigma_{C}(-x^{*})\ge-\langle x,x^{*}\rangle$ (see
\cite{MR1009594}). Hence
\begin{equation}
0\le\inf_{X\times X^{*}}F=-(\varphi_{T}+g)^{*}(0,0)=-\min_{(x,x^{*})\in X\times X^{*}}\{\psi_{T}(x,x^{*})+g^{*}(-x^{*},-x)\},\label{conv}
\end{equation}
because $g$ is continuous on $\operatorname*{int}C\times X^{*}$(see
f.i. \cite[Theorem 2.8.7, p. 126]{MR1921556}). Here $\psi_{T}(x,x^{*})=\varphi_{T}^{*}(x^{*},x)$,
$(x,x^{*})\in X\times X^{*}$; the convex conjugation being taken
with respect to the dual system $(X\times X^{*},X^{*}\times X^{**})$
and, for every $(x,x^{*})\in X\times X^{*}$, $\psi_{T}(x,x^{*})\ge\langle x,x^{*}\rangle$
because $T$ is monotone (see e.g. \cite[(12)]{MR2577332}).

From $g^{*}(x^{*},x)=\iota_{C}(-x)+\sigma_{C}(x^{*})$, $(x,x^{*})\in X\times X^{*}$
and (\ref{conv}) there exists $(\bar{x},\bar{x}^{*})\in X\times X^{*}$
such that $\psi_{T}(\bar{x},\bar{x}^{*})+\iota_{C}(\bar{x})+\sigma_{C}(-\bar{x}^{*})\le0$
which implies that $\iota_{C}(\bar{x})+\sigma_{C}(-\bar{x}^{*})=-\langle\bar{x},\bar{x}^{*}\rangle$,
i.e., $-\bar{x}^{*}\in N_{C}(\bar{x})$ and $\psi_{T}(\bar{x},\bar{x}^{*})=\langle\bar{x},\bar{x}^{*}\rangle$,
that is, $\bar{x}^{*}\in T(\bar{x})$ since $T$ is representable
(see \cite[Theorem 1, p. 270]{MR2577332}). Therefore $0\in(T+N_{C})(\bar{x},\bar{x}^{*})$
and so $0\in R(T+N_{C})$. \end{proof}

\strut

\begin{proof}[Proof of Theorem 1] First we prove that we can assume
without loss of generality that $D(B)$ is bounded. Indeed, assume that the result is true for that case.
Let $z=(x,x^{*})$ be m.r. to $A+B$. Take $C\subset X$ closed convex
and bounded with $x\in\operatorname*{int}C$ and $D(A)\cap\operatorname*{int}D(B)\cap\operatorname*{int}C\neq\emptyset$
e.g. $C:=[x_{0},x]+S$, where $[x_{0},x];=\{tx_{0}+(1-t)x\mid0\le t\le1\}$
and $S$ is a closed convex bounded neighborhood of $0$, and $x_{0}\in D(A)\cap\operatorname*{int}D(B)$.
Note that $z$ is m.r. to $A+B+N_{C}=A+(B+N_{C})$ which is
maximal monotone since, according to Theorem 2,  $B+N_{C}$
is maximal monotone,  $D(B+N_{C})$ is bounded, and
$x_{0}\in D(A)\cap\operatorname*{int}D(B+N_{C})\neq\emptyset$.
Hence $z\in\operatorname*{Graph}(A+B+N_{C})$ or $x^{*}\in(A+B)(x)$
because $N_{C}(x)=\{0\}$. Therefore $A+B$ is maximal monotone. 

It remains to prove that, whenever $D(B)$ is bounded,  $R(A+B)=X^{*}$ or
sufficiently $0\in R(A+B)$ since $A+B$ is representable (see again
\cite[Corollary\ 5.6]{MR2453098}). 

Let $F(x,x^{*}):=\varphi_{A}(x,x^{*})+\varphi_{B}(x,-x^{*})$, $g(x,x^{*}):=\varphi_{B}(x,-x^{*})$,
$(x,x^{*})\in X\times X^{*}$. Since $A,B$ are maximal monotone,
for every $(x,x^{*})\in X\times X^{*}$, $\varphi_{A}(x,x^{*}),\varphi_{B}(x,x^{*})\ge\langle x,x^{*}\rangle$
which imply $F\ge0$ and so
\[
0\le\inf_{X\times X^{*}}F=-(\varphi_{A}+g)^{*}(0,0)=-\min_{(x,x^{*})\in X\times X^{*}}\{\psi_{A}(x,x)+\psi_{B}(x,-x^{*})\},
\]
because $g$ is continuous on $\operatorname*{int}D(B)\times X^{*}$.

There exists $(\bar{x},\bar{x}^{*})\in X\times X^{*}$ such that $\psi_{A}(\bar{x},\bar{x}^{*})+\psi_{B}(\bar{x},-\bar{x}^{*})\le0$
which implies that $\psi_{A}(\bar{x},\bar{x}^{*})=\langle\bar{x},\bar{x}^{*}\rangle$,
$\psi_{B}(\bar{x},-\bar{x}^{*})=-\langle\bar{x},\bar{x}^{*}\rangle$,
i.e., $\bar{x}^{*}\in A(\bar{x})$ and $-\bar{x}^{*}\in B(\bar{x})$
from which $0\in R(A+B)$. \end{proof}

\begin{remark} Theorem 2 still holds if we replace the assumption $C$ bounded with $D(T)$ bounded. In this case an alternate proof of Theorem 1 can be performed with $A+N_C$ instead of $A$ and a similar argument as in the current proof. 
\end{remark}

\bibliographystyle{plain}

\end{document}